\documentclass[a4]{article}
\usepackage{amsmath}
\usepackage{amsfonts}
\usepackage{amssymb}
\usepackage{latexsym}
\usepackage{amsbsy}

\begin{document}

\newcommand{\nc}[2]{\newcommand{#1}{#2}}
\newcommand{\rnc}[2]{\renewcommand{#1}{#2}}
\rnc{\theequation}{\thesection.\arabic{equation}}
\def\note#1{{}}
\def\Label#1{\label{#1}\ifmmode\llap{[#1] }\else
\marginpar{\smash{\hbox{[#1]}}}\fi}

\newtheorem{definition}{Definition $\!\!$}[section]
\newtheorem{proposition}[definition]{Proposition $\!\!$}
\newtheorem{lemma}[definition]{Lemma $\!\!$}
\newtheorem{corollary}[definition]{Corollary $\!\!$}
\newtheorem{theorem}[definition]{Theorem $\!\!$}
\newtheorem{example}[definition]{\it Example $\!\!$}
\newtheorem{remark}[definition]{\it Remark $\!\!$}

\nc{\beq}{\begin{equation}}
\nc{\eeq}{\end{equation}}
\rnc{\[}{\beq}
\rnc{\]}{\eeq}
\rnc{\Label}{\label}
\nc{\ba}{\begin{array}}
\nc{\ea}{\end{array}}
\nc{\bea}{\begin{eqnarray}}
\nc{\beas}{\begin{eqnarray*}}
\nc{\eeas}{\end{eqnarray*}}
\nc{\eea}{\end{eqnarray}}
\nc{\be}{\begin{enumerate}}
\nc{\ee}{\end{enumerate}}
\nc{\bd}{\begin{diagram}}
\nc{\ed}{\end{diagram}}
\nc{\bi}{\begin{itemize}}
\nc{\ei}{\end{itemize}}
\nc{\bpr}{\begin{proposition}}
\nc{\bth}{\begin{theorem}}
\nc{\ble}{\begin{lemma}}
\nc{\bco}{\begin{corollary}}
\nc{\bre}{\begin{remark}}
\nc{\bex}{\begin{example}}
\nc{\bde}{\begin{definition}}
\nc{\ede}{\end{definition}}
\nc{\epr}{\end{proposition}}
\nc{\ethe}{\end{theorem}}
\nc{\ele}{\end{lemma}}
\nc{\eco}{\end{corollary}}
\nc{\ere}{\hfill\mbox{$\losenge$}\end{remark}}
\nc{\eex}{\hfill\mbox{$\losenge$}\end{example}}
\nc{\bpf}{{\it Proof.~~}}
\nc{\epf}{\hfill\mbox{$\square$}\vspace*{3mm}}
\nc{\hsp}{\hspace*}
\nc{\vsp}{\vspace*}

\nc{\ot}{\otimes}
\nc{\te}{\!\ot\!}
\nc{\bmlp}{\mbox{\boldmath$\left(\right.$}}
\nc{\bmrp}{\mbox{\boldmath$\left.\right)$}}
\nc{\LAblp}{\mbox{\LARGE\boldmath$($}}
\nc{\LAbrp}{\mbox{\LARGE\boldmath$)$}}
\nc{\Lblp}{\mbox{\Large\boldmath$($}}
\nc{\Lbrp}{\mbox{\Large\boldmath$)$}}
\nc{\lblp}{\mbox{\large\boldmath$($}}
\nc{\lbrp}{\mbox{\large\boldmath$)$}}
\nc{\blp}{\mbox{\boldmath$($}}
\nc{\brp}{\mbox{\boldmath$)$}}
\nc{\LAlp}{\mbox{\LARGE $($}}
\nc{\LArp}{\mbox{\LARGE $)$}}
\nc{\Llp}{\mbox{\Large $($}}
\nc{\Lrp}{\mbox{\Large $)$}}
\nc{\llp}{\mbox{\large $($}}
\nc{\lrp}{\mbox{\large $)$}}
\nc{\lbc}{\mbox{\Large\boldmath$,$}}
\nc{\lc}{\mbox{\Large$,$}}
\nc{\Lall}{\mbox{\Large$\forall\;$}}
\nc{\bc}{\mbox{\boldmath$,$}}
\nc{\ra}{\rightarrow}
\nc{\ci}{\circ}
\nc{\cc}{\!\ci\!}
\nc{\lra}{\longrightarrow}
\rnc{\to}{\mapsto}
\nc{\imp}{\Rightarrow}
\rnc{\iff}{\Leftrightarrow}
\nc{\inc}{\mbox{$\,\subseteq\;$}}
\rnc{\subset}{\inc}
\nc{\0}{\sp{(0)}}
\nc{\1}{\sp{(1)}}
\nc{\2}{\sp{(2)}}
\nc{\3}{\sp{(3)}}
\nc{\4}{\sp{(4)}}
\nc{\5}{\sp{(5)}}
\nc{\6}{\sp{(6)}}
\nc{\7}{\sp{(7)}}
\newcommand{\squareneu}{\square}
\def\tr{{\rm tr}}
\def\Tr{{\rm Tr}}
\def\st{\stackrel}
\def\<{\langle}
\def\>{\rangle}
\def\d{\mbox{$\mathop{\mbox{\rm d}}$}}
\def\b{\mbox{$\mathop{\mbox{\rm b}}$}}
\def\id{\mbox{$\mathop{\mbox{\rm id}}$}}
\def\ker{\mbox{$\mathop{\mbox{\rm Ker$\,$}}$}}
\def\hom{\mbox{$\mathop{\mbox{\rm Hom}}$}}
\def\corep{\mbox{$\mathop{\mbox{\rm Corep}}$}}
\def\im{\mbox{$\mathop{\mbox{\rm Im}}\,$}}
\def\map{\mbox{$\mathop{\mbox{\rm Map}}$}}
\def\o{\sp{[1]}}
\def\t{\sp{[2]}}
\def\mo{\sp{[-1]}}
\def\z{\sp{[0]}}

\nc{\spp}{\mbox{${\cal S}{\cal P}(P)$}}
\nc{\ob}{\mbox{$\Omega\sp{1}\! (\! B)$}}
\nc{\op}{\mbox{$\Omega\sp{1}\! (\! P)$}}
\nc{\oa}{\mbox{$\Omega\sp{1}\! (\! A)$}}
\nc{\dr}{\mbox{$\Delta_{R}$}}
\nc{\dsr}{\mbox{$\Delta_{\Omega^1P}$}}
\nc{\ad}{\mbox{$\mathop{\mbox{\rm Ad}}_R$}}
\nc{\as}{\mbox{$A(S^3\sb s)$}}
\nc{\bs}{\mbox{$A(S^2\sb s)$}}
\nc{\slc}{\mbox{$A(SL(2,\C))$}}
\nc{\suq}{\mbox{$\cO(SU_q(2))$}}
\nc{\tc}{\widetilde{can}}
\def\slq{\mbox{$\cO(SL_q(2))$}}
\def\asq{\mbox{$\cO(S_{q,s}^2)$}}
\def\esl{{\mbox{$E\sb{\frak s\frak l (2,{\Bbb C})}$}}}
\def\esu{{\mbox{$E\sb{\frak s\frak u(2)}$}}}
\def\ox{{\mbox{$\Omega\sp 1\sb{\frak M}X$}}}
\def\oxh{{\mbox{$\Omega\sp 1\sb{\frak M-hor}X$}}}
\def\oxs{{\mbox{$\Omega\sp 1\sb{\frak M-shor}X$}}}
\def\Fr{\mbox{Fr}}

\rnc{\epsilon}{\varepsilon}
\rnc{\phi}{\varphi}
\nc{\ha}{\mbox{$\alpha$}}
\nc{\hb}{\mbox{$\beta$}}
\nc{\hg}{\mbox{$\gamma$}}
\nc{\hd}{\mbox{$\delta$}}
\nc{\he}{\mbox{$\varepsilon$}}
\nc{\hz}{\mbox{$\zeta$}}
\nc{\hs}{\mbox{$\sigma$}}
\nc{\hk}{\mbox{$\kappa$}}
\nc{\hm}{\mbox{$\mu$}}
\nc{\hn}{\mbox{$\nu$}}
\nc{\hl}{\mbox{$\lambda$}}
\nc{\hG}{\mbox{$\Gamma$}}
\nc{\hD}{\mbox{$\Delta$}}
\nc{\hT}{\mbox{$\Theta$}}
\nc{\ho}{\mbox{$\omega$}}
\nc{\hO}{\mbox{$\Omega$}}
\nc{\hp}{\mbox{$\pi$}}
\nc{\hP}{\mbox{$\Pi$}}

\nc{\qpb}{quantum principal bundle}
\def\gal{-Galois extension}
\def\hge{Hopf-Galois extension}
\def\ses{short exact sequence}
\def\csa{$C^*$-algebra}
\def\ncg{noncommutative geometry}
\def\wrt{with respect to}
\def\Ha{Hopf algebra}
\newcommand{\ayd}{anti-Yetter-Drinfeld}

\def\C{{\Bbb C}}
\def\N{{\Bbb N}}
\def\R{{\Bbb R}}
\def\Z{{\Bbb Z}}
\def\T{{\Bbb T}}
\def\Q{{\Bbb Q}}
\def\cO{{\mathcal O}}
\def\O{\cO}
\def\cT{{\cal T}}
\def\cK{{\cal K}}
\def\cH{{\cal H}}
\def\ch{{\cal H}}
\def\H{{\cal H}}
\def\ta{\tilde a}
\def\tb{\tilde b}
\def\td{\tilde d}
\newcommand{\weg}[1]{}

{\noindent\LARGE\bf
Hopf-cyclic homology and cohomology with\\ coefficients }\\
\medskip

\noindent{\large\bf
Piotr M.~Hajac$^a$, Masoud Khalkhali$^b$,
 Bahram Rangipour$^b$,\\
 Yorck Sommerh\"auser$^c$\\
\medskip\medskip

\footnotesize\bf\begin{itemize}
\vspace*{-4mm}\item[$^a$]
Instytut Matematyczny, Polska Akademia Nauk,
ul.\ \'Sniadeckich 8, Warszawa, 00-956 Poland,
and
Katedra Metod Matematycznych Fizyki, Uniwersytet Warszawski,
ul.\ Ho\.za 74, Warszawa, 00-682 Poland;
URL: http://www.fuw.edu.pl/$\!\widetilde{\phantom{m}}\!$pmh\\

\vspace*{-4mm}\item[$^b$]
Department of Mathematics, University of Western Ontario, London ON, Canada;
E-mail: masoud@uwo.ca, brangipo@uwo.ca\\

\vspace*{-4mm}\item[$^c$]
Mathematisches Institut, Universit\"at M\"unchen,
Theresienstr.\ 39, 80333 M\"unchen,  Germany;
E-mail: sommerh@mathematik.uni-muenchen.de
\end{itemize}
}
\bigskip

\ $\!\!$\hrulefill\
\medskip\smallskip

\noindent{\footnotesize
\vspace*{-.5mm}
{\bf Abstract.}
Following the idea of an invariant differential complex,
we construct
general-type cyclic modules that provide
\vspace*{-.5mm}
 the common denominator of known cyclic
theories. The cyclicity of these modules is governed by Hopf-algebraic
\vspace*{-.5mm}
 structures.
We prove that the existence of a cyclic operator forces a modification of the
Yetter-Drinfeld compatibility
\vspace*{-.5mm}
 condition leading to the concept of a stable
anti-Yetter-Drinfeld module.
This module plays the role of the space of
\vspace*{-.5mm}
coefficients in the thus obtained
cyclic cohomology of module algebras and coalgebras, and
the cyclic homology and
\vspace*{-.5mm}
cohomology of comodule algebras. Along the lines of Connes and Moscovici,
we show that there is a pairing between
\vspace*{-.5mm}
the cyclic cohomology of a module
coalgebra acting on a module algebra and closed  0-cocycles on the latter.
The
\vspace*{-.5mm}
pairing takes values in the
usual cyclic cohomology of the algebra. Similarly, we argue that
there is an analogous
\vspace*{-.5mm}
pairing between closed  0-cocycles
of a module coalgebra and the cyclic cohomology of a module algebra.
}\ ~\\

{\noindent\normalsize\bf
Homologie et cohomologie Hopf-cycliques \`a coefficients}
\medskip

\noindent{\footnotesize
\vspace*{-.5mm}
{\bf R\'esum\'e.}
Suivant l'id\'ee d'un complexe diff\'erentiel invariant, nous
construisons des modules cycliques de type
\vspace*{-.5mm}
 g\'en\'eral qui fournissent
un d\'enominateur commun aux th\'eories cycliques connues. Le
caract\`ere cyclique de ces
\vspace*{-.5mm}
modules est gouvern\'e par des structures
Hopf-alg\'ebriques. Nous montrons que l'existence d'un op\'erateur
cyclique
\vspace*{-.5mm}
oblige \`a une modification de la condition de
compatibilit\'e de Yetter-Drinfeld et m\`ene au concept de module
anti-
\vspace*{-.5mm}
Yetter-Drinfeld stable. Ce module joue le r\^ole d'espace de
coefficients pour la cohomologie de modules alg\`ebres  et
\vspace*{-.5mm}
de modules cog\`ebres
 ainsi obtenue, ainsi que pour l'homologie et la cohomologie
cycliques de comodules alg\`ebres.
\vspace*{-.5mm}
Comme l'ont fait Connes et Moscovici
pour leur cohomologie, nous montrons qu'il existe un appariement entre
la
\vspace*{-.5mm}
cohomologie cyclique d'un module cog\`ebre agissant sur un module
alg\`ebre et les 0-cycles ferm\'es sur ce dernier. Cet
\vspace*{-.5mm}
appariement
prend ses valeurs dans la cohomologie cyclique usuelle de
l'alg\`ebre. De fa\c con similaire, nous \'etablissons
\vspace*{-.5mm}
un appariement
analogue entre les 0-cycles ferm\'es d'un module cog\`ebre et la
cohomologie cyclique d'un module
\vspace*{-.5mm}
alg\`ebre.
}

\ $\!\!$\hrulefill\
\bigskip

\noindent{\bf Version fran\c{c}aise abr\'eg\'ee.}

Soit $H$ une alg\`ebre de Hopf avec une antipode bijective, et soit
$M$ \`a la fois un module et un comodule sur $H.$ Nous disons que $M$
est un module {\em \ayd} s'il satisfait \`a l'une des conditions \'equivalentes
(\ref{ayd1})--(\ref{ayd4}). Nous disons que $M$ est {\em stable} s'il
a la propri\'et\'e  $action\ci coaction=\id$. Le r\'esultat principal
de cet article est la construction de modules cycliques qui nous
permettent de d\'efinir la cohomologie cyclique de $H$-modules
cog\`ebres et de $H$-modules alg\`ebres \`a coefficients dans des modules
{\ayd} stables. De plus, si un $H$-module cog\`ebre $C$ agit sur
un $H$-module alg\`ebre $A$ au sens o\`u nous avons une application
lin\'eaire $C\ot A\ra A$ satisfaisant $c(ab)=(c^{(1)}a)(c^{(2)}b)$,
$c1=\he(c)$, $h(ca)=(hc)a$, pour tous $h\in H$, $c\in C$, alors il
existe un certain appariement entre ces cohomologies pour $C$ et $A$.
\bth
Soit $H$ une alg\`ebre de Hopf avec une antipode bijective et soit $M$
un module {\ayd} sur $H$ stable. Alors les formules (\ref{f1})--(\ref{tau}) (resp.\
(\ref{a})--(\ref{d})) d\'efinissent une structure de  module cyclique
sur $\{M\ot_HC^{\ot(n+1)}\}_{n\in\N}$ (resp.\
$\{\hom^H(M\ot A^{\ot(n+1)},k)\}_{n\in\N}$)
 pour tout $H$-module cog\`ebre $C$ (resp. alg\`ebre $A$). Si de plus
 $C$ agit sur $A$ et  $\Q\inc k$, alors, pour tout  $n\in\N$, les
 formules  (\ref{pair1}) et (\ref{pair2}) d\'efinissent,
 respectivement, les homomorphismes suivants:
\beq
HC^n_H(C,M)\ot HC^0_H(A,M)\st{\#_{n,0}}{\lra} HC^n(A),~~~
HC^0_H(C,M)\ot HC^n_H(A,M)\st{\#_{0,n}}{\lra} HC^n(A).
\eeq
\ethe\noindent
De fa\c con similaire, pour tout comodule alg\`ebre, les formules
(\ref{coh1})--(\ref{coh4})  et (\ref{hom1})--(\ref{hom4})
d\'eterminent respectivement notre construction d'une homologie et
d'une cohomologie cycliques \`a coefficients dans un module {\ayd}
stable.

\section*{Introduction}
\setcounter{equation}{0}

Ever since its invention, among main applications of cyclic cohomology was
the computation of $K$-theoretical invariants. While enhancing the feasibility
of such computations, Connes and Moscovici discovered a new type of cyclic
cohomology, notably the cyclic cohomology of Hopf algebras \cite{cm98}.
Invariant cyclic homology  introduced in \cite{kr} generalizes the
Connes-Moscovici theory
and its dual  version \cite{kr02}. It shows that  passage from the cyclic homology
of algebras to the cyclic cohomology of Hopf algebras is remarkably  similar to
passage from de Rham cohomology  to cohomology of Lie algebras via invariant
de Rham cohomology \cite{ce48}. The idea of employing invariant complexes
proved to be a key in resolving by significantly more effective means the
technical challenge
of showing that the $(n+1)$-power of the cyclic operator $\tau_n$ is identity
\cite[p.102]{cm99}, and allowed the introduction of higher-dimensional coefficients.

We continue this strategy herein. Our motivation is to obtain and understand
computable invariants of $K$-theory. The aim of this paper is to provide
general framework for cyclic theories whose cyclicity is based on Hopf-algebraic
structures. We refer to such homology and cohomology as {\em Hopf-cyclic}.
First we define and provide sources of examples
of stable \ayd\ modules that play the role of coefficients for Hopf-cyclic theory.
In particular, we claim that modular pairs in involution are precisely 1-dimensional
stable \ayd\ modules.
Then we construct cyclic module structures on invariant complexes for module
coalgebras and module algebras, respectively.
It turns out that  the cyclic cohomology of Hopf algebras
is a special case of the former, whereas both twisted \cite{kmt03} and usual cyclic
cohomology are special cases of the latter. As a bonus of generality, we obtain
now a very short proof of Connes-Moscovici key result \cite[Theorem~1]{cm99}.
Furthermore, as \hd-invariant
\hs-traces can be viewed as closed 0-cocycles on a module algebra, our pairing
for Hopf-cyclic cohomology generalizes the Connes-Moscovici
transfer map \cite[Proposition~1]{cm99} from the cyclic cohomology of Hopf algebras
to ordinary cyclic cohomology. Finally, we end this article by deriving Hopf-cyclic
homology and cohomology of comodule algebras. This extends the formalism
for comodule algebras provided in \cite{kr}.

Partly for the sake of simplicity, throughout the paper we assume that $H$
is a Hopf algebra with a bijective antipode.
 On one hand side, the bijectivity of the antipode is implied
by the existence of a modular pair in involution, so that then it need not be
assumed. On the other hand, some parts of arguments might work even if the
antipode is not bijective. We avoid such discussions.
The  coproduct, counit and antipode of $H$ are denoted by
$\Delta$,  $\epsilon$ and $S$,  respectively.
For the coproduct we use the  notation $\Delta(h)=h^{(1)}\ot h^{(2)}$,
for a left coaction on $M$ we write $_M\Delta(m)=m^{(-1)}\ot m^{(0)}$, and for
a right coaction  $\Delta_M(m)=m^{(0)}\ot m^{(1)}$.  The summation
 symbol is suppressed everywhere. We assume all algebras to be associative, unital
and over the same ground field $k$. The symbol ${\cal O}(X)$ stands for the
algebra of polynomial functions  on $X$.

\section{Stable anti-Yetter-Drinfeld modules}
\setcounter{equation}{0}

It turns out that, in order to incorporate coefficients into  cyclic theory,
we  need to alter the concept of a Yetter-Drinfeld module by replacing
  the antipode by its inverse in the Yetter-Drinfeld compatibility condition
 between actions and coactions. We call the modules-comodules satisfying
 the thus modified  Yetter-Drinfeld compatibility condition {\em\ayd\
 modules}\footnote{
 This concept was  devised independently by Ch.~Voigt
and, also independently, by P.~Jara and D.~Stefan.}.
   Just as Yetter-Drinfeld modules
 come in 4 different versions depending on the side of actions and coactions
 (see \cite[p.181]{cmz02} for a general formulation),
 so do the \ayd\ modules.
All versions are completely equivalent and can be derived from one another
 by replacing a Hopf algebra
 $H$ by $H^{cop}$, $H^{op}$, or $H^{op,cop}$, respectively.
 Explicitly, we have the following:
\begin{definition}\Label{ayd}
 Let $H$ be a Hopf algebra with a bijective antipode $S$,
 and $M$ a module and comodule
 over $H$. We call $M$ an \ayd\ module iff
 the action and coaction are compatible in the following sense:
\bea
&&\label{ayd1}
_M\Delta(hm)=h^{(1)}m^{(-1)}S^{-1}(h^{(3)})\ot h^{(2)}m^{(0)}
~~~
\mbox{\em if $M$ is a left module and a left comodule},
~~~~~~\\ &&
\Delta_M(hm)=h^{(2)}m^{(0)}\ot h^{(3)}m^{(1)}S(h^{(1)})
~~~
\mbox{\em if $M$ is a left module and a right comodule},
~~~~~~\\ &&
_M\Delta(mh)=S(h^{(3)})m^{(-1)}h^{(1)}\ot m^{(0)}h^{(2)}
~~~
\mbox{\em if $M$ is a right module and a left comodule},
~~~~~~\\ &&\label{ayd4}
\Delta_M(mh)= m^{(0)}h^{(2)}\ot S^{-1}(h^{(1)})m^{(1)}h^{(3)}
~~~
\mbox{\em if $M$ is a right module and a right comodule}.
~~~~~~
\eea
\end{definition}\noindent
To make  cyclic theory work, we also need to assume that the action splits
coaction, i.e., for all $m\in M$,
$m^{(-1)}m^{(0)}=m$, $m^{(1)}m^{(0)}=m$, $m^{(0)}m^{(-1)}=m$,
$m^{(0)}m^{(1)}=m$, for the left-left, left-right, right-left, and right-right
version respectively. We call modules satisfying this condition {\em stable}.
Let us emphasize that it is the anti-Yetter-Drinfeld condition rather than
the Yetter-Drinfeld condition that makes the homomorphism
$action\ci coaction$ $H$-linear and $H$-colinear. Therefore the stability condition
$action\ci coaction=\id$ suits the former and not the latter.
The first class of examples of stable \ayd\ modules is provided by modular
pairs in involution \cite[p.8]{cm00}.
\ble
Let the ground field
$k$ be a right module over $H$ via a character $\hd$ and a left comodule
over $H$ via a group-like $\hs$. Then $k=^\sigma\!\!\!k_\delta$ is a stable \ayd\ module
{\em if and only if} $(\hd,\hs)$ is a modular pair in involution.
\ele\noindent
This combined with the following lemma guarantees a rich source of examples of
\ayd\ modules.
\ble
Let $N$ be a Yetter-Drinfeld module and $M$ an \ayd~ module.
Then $N\otimes M$  is an \ayd~ module via
$h(n\ot m)=h^{(1)}n\ot h^{(2)}m$,
$_{N\ot M}\hD(n\ot m)=n^{(-1)}m^{(-1)}\ot n^{(0)}\ot m^{(0)}$,
for the left-left case, and via
$h(n\ot m)=h^{(2)}n\ot h^{(1)}m$,
$\hD_{N\ot M}(n\ot m)=n^{(0)}\ot m^{(0)}\ot n^{(1)}m^{(1)}$,
for the left-right case. Similarly, $M\ot N$ is an \ayd\ module via
$(m\ot n)h=mh^{(2)}\ot nh^{(1)}$,
$_{M\ot N}\hD(n\ot m)=m^{(-1)}n^{(-1)}\ot m^{(0)}\ot n^{(0)}$,
for the right-left case, and via
$(m\ot n)h=mh^{(1)}\ot nh^{(2)}$,
$\hD_{M\ot N}(m\ot n)=m^{(0)}\ot n^{(0)}\ot m^{(1)}n^{(1)}$,
for the left-right case.
\ele
An intermediate step between modular pairs in involution and stable \ayd\ modules
is given by matched and comatched pairs of \cite{kr}.
Whenever the antipode is equal to its inverse, the difference between the
Yetter-Drinfeld and \ayd\ conditions disappears.
For a group ring Hopf algebra $kG$, a left
$H$-comodule
 is simply a $G$-graded vector space $M={\bigoplus}_{g\in G}M_g$,
  where the coaction is defined by $M_g\ni m\mapsto g\otimes m$.
 An action of $G$ on $M$ defines an \mbox{(anti-)}
Yetter-Drinfeld module if and only if
for all $g,h\in G$ and $m\in M_g$ we have
 $hm\in M_{hgh^{-1}}$. The stability condition means simply that $gm=m$ for all
$g\in G$, $m\in M_g$.
A very concrete classical example of a stable (anti-)Yetter-Drinfeld module is
provided by the Hopf fibration. Then $H=\cO(SU(2))$ and $M=\cO(S^2)$. Since
$S^2\cong SU(2)/U(1)$, we have a natural left action of $SU(2)$ on $S^2$.
Its pull-back
makes $M$ a left $H$-comodule. On the other hand, one can view $S^2$ as the set
of all traceless matrices of $SU(2)$. The pull-back of this embedding
$j:S^2{\hookrightarrow}SU(2)$ together with the multiplication in $\cO(S^2)$
defines a left action of $H$ on $M$. It turn out that the equivariance property
$j(gx)=gj(x)g^{-1}$ guarantees the \ayd\ condition, and this combined with
the injectivity of $j$ ensures the stability of $M$. This stability mechanism
can be generalized in the following way.
\ble
 Let $M$ be an algebra and a left $H$-comodule. Assume that
   $\pi: H\ra M$ is an epimorphism of algebras and
   the action $hm=\pi(h)m$ makes $M$ an \ayd\ module. Assume
also that $\pi(1^{(-1)})1^{(0)}=1$. Then
   $M$ is a stable module.
\ele

Another source of examples is provided by Hopf-Galois theory. These examples are
purely quantum in the sense that the employed actions are automatically trivial
for commutative algebras. To fix the notation and terminology, recall that an algebra
and an $H$-comodule is called a comodule algebra if the coaction is an algebra
homomorphism. An $H$-extension $B:=\{p\in P\;|\;\hD_P(p)=p\ot 1\}\inc P$ is called
Hopf-Galois iff the canonical map $can: P\ot_BP\ra P\ot H$, $can(p\ot p')=
p\hD(p')$, is bijective.
The bijectivity assumption
allows us to define the translation map $T:H\ra P\ot_BP$, $T(h):=can^{-1}(1\ot h)
=:h^{[1]}\ot_Bh^{[2]}$ (summation suppressed). It can be shown that when
everything is over a field (our standing assumption), the centralizer
$Z_B(P):=\{p\in P\;|\;bp=pb, \forall b\in B\}$ of $B$ in $P$ is a subcomodule of $P$.
On the other hand,
 the formula $ph=h^{[1]}ph^{[2]}$ defines a right action on $Z_B(P)$ called
the Miyashita-Ulbrich action. This action and coaction satisfy the
Yetter-Drinfeld compatibility condition \cite[(3.11)]{dt89}. The following proposition modifies
the Miyashita-Ulbrich action  so as to obtain stable \ayd\ modules.
\bpr\Label{mu}
Let  $B\inc P$ be a Hopf-Galois
$H$-extension such that $B$ is central in $P$. Then $P$ is a right-right
stable \ayd\ module via the action
$
ph=(S^{-1}(h))^{[2]}p(S^{-1}(h))^{[1]}
$
and the right coaction on $P$.
\epr\noindent
The simplest examples are obtained for $P=H$. A broader class is given by the
so-called Galois objects \cite{c-s98}. Then quantum-group coverings at roots of unity
provide examples with central coinvariants bigger than the ground field
(see \cite{dhs99} and examples therein).

\section{Cyclic cohomology of module algebras and coalgebras}
\setcounter{equation}{0}

An algebra $A$ which is a module over a Hopf algebra $H$ such that
$h(ab)=(h^{(1)}a)(h^{(2)}b)$ and $h1=\he(h)$, for all $h\in H$, $a,b\in A$,
is called an {\em $H$-module algebra}. Similarly, a coalgebra $C$
which is a module over a Hopf algebra $H$ such that
$\hD(hc)=(h^{(1)}c^{(1)})\ot(h^{(2)}c^{(2)})$ and
$\he(hc)=\he(h)\he(c)$, for all $h\in H$ and $c\in C$,
is called an {\em $H$-module coalgebra}.
In this section we construct cyclic modules for both module algebras and coalgebras.

We begin with the coalgebra case. First, we take a left $H$-comodule $M$,
$H$-module coalgebra $C$, and
 define an auxiliary simplicial complex
${\cal C}^n(C,M):=M\ot C^{\ot(n+1)}=(M\ot C)\ot C^{\ot n}$, $n\in\N$,
by endowing $M\ot C$ with the bicomodule structure
$_{M\ot C}\hD(m\ot c)=m^{(-1)}c^{(1)}\ot(m^{(0)}\ot c^{(2)})$,
$\hD_{M\ot C}(m\ot c)=(m\ot c^{(1)})\ot c^{(2)}$. Then the standard formulas for
faces and degeneracies on the complex $\{\mbox{bicomodule}\ot C^{\ot n}\}_{n\in\N}$
 translate immediately into
\begin{eqnarray}\Label{f1}
\delta_i(m\ot c_0 \otimes \dots \otimes c_{n-1})
\!\!\!&=&\!\!\!
m\ot  c_0 \otimes\dots
\otimes  c_i^{(1)}\otimes c_i^{(2)}\otimes c_{n-1},
 ~~~ 0 \leq i <n ~~~\mbox{(faces)},~~~~~~~~~
\\
\delta_{n}(m\ot c_0 \otimes \dots \otimes c_{n-1})
\!\!\!&=&\!\!\!
m^{(0)}\ot  c_0^{(2)}\otimes c_1
\otimes \dots \otimes c_{n-1}  \otimes m^{({-1})}c_0^{(1)}
~~~ \mbox{(flip-over face)},~~~~~~~~~
\\
\sigma_i(m\ot c_0 \otimes \dots \otimes c_{n+1})
\!\!\!&=&\!\!\!
m\ot c_0 \otimes \dots
 \otimes \varepsilon(c_{i+1})\otimes
\dots\otimes c_{n+1},~~~0\leq i \leq n~~~\mbox{(degeneracies)}.~~~~~~
\end{eqnarray}
It is  straightforward to check that these morphisms together with
\beq\Label{tau}
\tau_n(m\ot c_0 \otimes  \dots \otimes c_n)
=
m^{(0)}\ot c_1 \otimes
\dots \otimes c_n \otimes m^{({-1})}c_0
\eeq
form a paracyclic structure on $\{{\cal C}^n(C,M)\}_{n\in\N}$.
Now let us assume that $M$ is also a right $H$-module. We can treat $C^{\ot(n+1)}$
as a left $H$-module via the diagonal action
(i.e.,
$h(c_0 \otimes  \dots \otimes c_n)=h^{(1)}c_0 \otimes  \dots \otimes h^{(n+1)}c_n$)
and define the quotient (invariant) complex ${\cal C}^n_H(C,M):=M\ot_H C^{\ot(n+1)}$,
$n\in\N$. Except for $\tau_n$ and $\hd_n$ it is clear that the aforementioned morphisms
are well defined on the quotient complex. The key result of this paper is that
$\tau_n$ is well defined if and only if $M$ is an \ayd\ module.
 More precisely, we have:
\bth
Let
$M$ be a left $H$-comodule and
a right $H$-module. Then
the map $\tau_n$ given by the formula (\ref{tau}) is well defined on
$M\ot_H C^{\ot(n+1)}$ for any $n\in\N$ and any $H$-module
coalgebra $C$ {\em if and only if}
$M$ is an \ayd\ module. If furthermore $M$ is stable, then
$\{{\cal C}^n_H(C,M)\}_{n\in\N}$ is a cyclic
module via (\ref{f1})--(\ref{tau}).
\ethe
\bpf
Note first that $\tau_n$ is well defined on $M\ot_H C^{\ot(n+1)}$ if and only if
\beq\Label{1}
m^{(0)}\ot_H (h^{(2)}(c_1\ot\dots\ot c_n)\ot m^{(-1)}h^{(1)}c_0)
=
(hm)^{(0)}\ot_H (c_1\ot\dots\ot c_n\ot (hm)^{(-1)}c_0).
\eeq
The equality $m^{(0)}\ot_H (h^{(2)}\ot m^{(-1)}h^{(1)})
=
(hm)^{(0)}\ot_H (1\ot (hm)^{(-1)})$ evidently implies (\ref{1}) for any $n$ and $C$,
and follows if we assume that  (\ref{1}) holds for any $n$ and $C$. (E.g., take
$n=2$, $C=H$ and $c_0=1=c_1$.) Now, observe that there is an isomorphism
$\Phi:H\ot H\ra H\ot H$, $\Phi(h\ot k)=h^{(1)}\ot S(h^{(2)})k$,
$\Phi^{-1}(h\ot k)=h^{(1)}\ot h^{(2)}k$. $\Phi$ is left $H$-linear if we view
the domain as a left $H$-module via the diagonal action and the counter-domain
as a left $H$-module via the multiplication in the first factor. Thus we have
an isomorphism
$M\ot_H(H\ot H)
\st{\mbox{\scriptsize id}\ot_H\Phi}{\lra}
(M\ot_HH)\ot H\cong M\ot H$.
Applying this isomorphism to the equality below the equation (\ref{1}) yields
$m^{(0)}h^{(2)}\ot S(h^{(3)})m^{(-1)}h^{(1)}
=(hm)^{(0)}\ot (hm)^{(-1)}$. This is equivalent to the \ayd\ condition, as desired.
Next, since $\hd_n=\tau_n\hd_0$, all morphisms are well defined on ${\cal C}^n_H(C,M)$,
if $M$ is an \ayd\ module. If $M$ is also stable, then the equality
$
(\tau_n)^{n+1}(m\ot_H (c_0 \otimes \dots \otimes c_{n}))
=m^{(0)}m^{(-1)}\ot_H (c_0 \otimes \dots \otimes c_{n})
$
entails that $(\tau_n)^{n+1}=\id$, as needed.
\epf\\
\noindent
For $C=H$ and $M=^{\sigma}\!\!\!k_{\delta}$, the complex
 $\{{\cal C}^n_H(C,M)\}_{n\in\N}$
becomes the cyclic module of Connes-Moscovici \cite{cm99}. As an intermediate
step, one can take the Hopf cotriples of \cite{kr}.

We can proceed much the same way in the algebra case. Again, we first have an
auxiliary complex ${\cal C}^n(A,M):=\hom(M\ot A^{\ot(n+1)},k)\cong
\hom(A^{\ot n},\hom(M\ot A,k))$. Here $A$ is a left
\linebreak
 $H$-module algebra and $M$
a left $H$-comodule. The formulas $b(m\ot a)=m^{(0)}\ot (S^{-1}(m^{(-1)})b)a$,
$(m\ot a)b=m\ot ab$ turn $M\ot A$, and consequently $\hom(M\ot A,k)$, into a bimodule.
Thus we can use the standard formulas \cite[p.37]{l-jl}
to define a simplicial
structure on $\{{\cal C}^n(A,M)\}_{n\in\N}$. To define an invariant subcomplex,
we assume now $M$ to be also a right $H$-module,
view $A^{\ot(n+1)}$ as a left $H$-module via the diagonal action,
$M\ot A^{\ot(n+1)}$ as a right $H$-module via
$(m\ot\widetilde{a})h=mh^{(1)}\ot S(h^{(2)})\widetilde{a}$, and $k$ as
a right $H$-module via the counit map. Then we can put
${\cal C}^n_H(A,M):=\hom_H(M\ot A^{\ot(n+1)},k)$, $n\in\N$. Also along the same
lines, one can prove that if $M$ satisfies the \ayd\ condition, then the morphisms
\bea
&&\Label{a}
(\delta_if)(m\ot a_0\ot \dots\ot a_n)
= f(m\ot a_0\ot \dots \ot a_i a_{i+1}\ot\dots \ot a_n),
~~~ 0 \leq i <n,~~~~~~
\\ && \Label{b}
(\delta_nf)(m\ot a_0\ot \dots\ot a_n)
=f(m^{(0)}\ot (S^{-1}(m^{(-1)})a_n)a_0\ot a_1\ot\dots \ot a_{n-1}),~~~~~~
\\ &&\Label{c}
(\sigma_if)(m\ot a_0\ot \dots\ot a_n)
= f(m\ot a_0 \ot \dots \ot a_i\ot 1 \ot \dots \ot a_n),
~~~0\le i\le n,~~~~~~
\\  &&\Label{d}
(\tau_nf)(m\ot a_0\ot \dots\ot a_n)
= f(m^{(0)}\ot (S^{-1}(m^{(-1)})a_n)\ot a_0\ot \dots \ot a_{n-1}),~~~~~~
\eea
define a paracyclic  structure on $\{{\cal C}^n_H(A,M)\}_{n\in\N}$. Adding
the stability assumption on $M$, one obtains:
\bth
If $M$ is a stable \ayd\ module and $A$ is a left $H$-module algebra, then
the maps (\ref{a})-(\ref{d}) define a cyclic module structure on
$\{{\cal C}^n_H(A,M)\}_{n\in\N}$.
\ethe\noindent
For $H=k=M$ we recover the usual cyclic cohomology of algebras. For
the Laurent polynomial Hopf algebra ($H=k[\sigma,\sigma^{-1}]$) and
$M=^\sigma\!\!\!k_{\varepsilon}$ we obtain the twisted cyclic cohomology
\cite{kmt03}. We can also take as a stable \ayd\ module a Hopf algebra $K$
thought of as a left comodule over itself via the coproduct, and as a right
module over itself via $k\cdot h=S(h^{(2)})kh^{(1)}$.  Then we arrive at
the construction considered in \cite{ak}. This is a special case of the construction
in Proposition~\ref{mu}: $P=K^{cop}=H$.

Finally, if we take $M=^\sigma\!\!\!k_{\delta}$,
 we can view a $\hd$-invariant
$\sigma$-trace  of \cite[Definition~1]{cm99} as a closed 0-cocycle in our complex.
On the other hand, if $\Q\inc k$, we can define
 cyclic cohomology as follows: $\b_n:=\sum_{i=0}^n(-1)^i\hd_i$,
${\cal Z}^n_H(*,M):=\ker \b_n|_{\mbox{\scriptsize\rm Ker}(\tau_n-(-1)^n)}$,
${\cal B}^n_H(*,M):=\im \b_{n-1}|_{\mbox{\scriptsize\rm Ker}(\tau_{n-1}-(-1)^{n-1})}$,
$HC^n_H(*,M):={\cal Z}^n_H(*,M)/{\cal B}^n_H(*,M)$, $n\in\N$.
This brings us
to the following generalization of  \cite[Proposition~1]{cm99}:
\bpr\Label{pair}
Let $C$ be a left $H$-module coalgebra, $A$ a left $H$-module algebra,
and $M$ a stable right-left \ayd\ over $H$. Assume also that $\Q\inc k$ and
$C$ acts on $A$ in
the sense that we have a linear map $C\ot A\ra A$ satisfying
$c(ab)=(c^{(1)}a)(c^{(2)}b)$, $c1=\he(c)$, $h(ca)=(hc)a$, for all $h\in H$, $c\in C$,
$a,b\in A$.
Then, for all $n\in\N$, the formulas
\beq\label{pair1}
([m\ot_H(c_0\ot\dots\ot c_n)]\# f)(a_0\ot\dots\ot a_n)
=f(m\ot_H(c_0a_0)\dots (c_na_n))
\eeq
\beq\label{pair2}
((m\ot_Hc)\#[f])(a_0\ot\dots\ot a_n)
= f(m\ot_H(c^{(1)}a_0)\ot\dots\ot c^{(n+1)}a_{n}))
\eeq
define, respectively, the following  homomorphisms:
\beq
HC^n_H(C,M)\ot {\cal Z}^0_H(A,M)\st{\#_{n,0}}{\lra} HC^n(A),~~~
{\cal Z}^0_H(C,M)\ot HC^n_H(A,M)\st{\#_{0,n}}{\lra} HC^n(A).
\eeq
\epr\noindent
We conjecture that one can construct along these line a cup product
between the Hopf-cyclic cohomology of a module coalgebra and
the Hopf-cyclic cohomology of a module algebra. This would
generalize in some sense the cup product in \cite[p.105]{c-a85}.

\section{Cyclic theory for comodule algebras}
\setcounter{equation}{0}

 In this section we  define  cyclic cohomology  with coefficients
  in a right-right  stable \ayd\ module, and cyclic homology  with coefficients
  in a left-left stable \ayd\ module.
The latter case generalizes  theory introduced \cite{kr}.
In both cases, we assume that $A$ is a right
$H$-comodule algebra.
Our strategy for constructing  cyclic theory is much as before:
bimodule structure $\leadsto$ simplicial structure  $\leadsto$ paracyclic structure
$\leadsto$ invariant complex $\leadsto$ cyclic theory. In the cohomology case,
we define the invariant subcomplex as ${\cal C}^{n,H}(A,M):=\hom^H(A^{\ot(n+1)},M)$,
$n\in\N$. Here   $M$ is a right-right stable \ayd\   module, and $A^{\ot(n+1)}$ is
viewed as a right comodule via the diagonal coaction.
 The following maps define, respectively, faces, degeneracies and cyclic
   operators on $\{{\cal C}^{n,H}(A,M)\}_{n\in \N}$.
\bea
&&\Label{coh1}
(\delta_if)(a_0\ot\dots\ot a_n)= f(a_0\ot \dots
 \ot a_i a_{i+1}\ot\dots \ot a_{n}),~~~
0\le i< n,~~~~~~
\\ &&
(\delta_{n}f)(a_0\ot\dots\ot a_n)
= f(a_{n}^{(0)}a_0\ot a_1\ot \dots
 \ot a_{n-1})a_{n}^{(1)},
\\  &&
(\sigma_if)(a_0\ot\dots\ot a_{n})= f(a_0 \ot \dots \ot a_i\ot
 1 \ot \dots \ot a_{n}),~~~
0\le i< n,~~~~~~
\\   &&\Label{coh4}
(\tau_n f)(a_0\ot\dots\ot a_n)= f(a_n^{(0)}\ot
  a_0\ot\dots \ot a_{n-1})a_n^{(1)}.
\eea
In the homology case, we take a left-left stable \ayd\ module $M$ and define the
invariant subcomplex via the cotensor product:
${\cal C}_{n}^H(A,M):=A^{\ot(n+1)}\square_HM$,
$n\in\N$.
(Recall that $X\square_HY:=\ker(\hD_X\ot\id-\id\ot\,_Y\hD)$.)
Here,   again, we view  $A^{\ot(n+1)}$ as a right $H$-comodule via
the diagonal coaction. The following homomorphisms establish, respectively,
the desired faces, degeneracies and cyclic operators on
$\{{\cal C}_{n}^H(A,M)\}_{n\in \N}$.
\bea
&&\Label{hom1}
d_i(a_0\ot\dots\ot a_n\ot m)= a_0\ot \dots
 \ot a_i a_{i+1}\ot\dots \ot a_n\ot m,~~~
0\le i< n,~~~ ~~~
\\ &&
d_n(a_0\ot\dots\ot a_n\ot m)= a_n^{(0)}a_0\ot a_1\dots
\ot a_{n-1}\ot a_{n}^{(1)}m,
\\  &&
s_i( a_0\ot\dots\ot a_n \ot m)= a_0 \ot \dots \ot a_i\ot
 1 \ot \dots \ot a_n\ot m,~~~
0\le i\le n,~~~ ~~~
\\ &&\Label{hom4}
t_n(a_0\ot\dots\ot a_n\ot m)= a_n^{(0)}\ot
  a_0\ot\dots \ot a_{n-1}\ot a_n^{(1)}m.
\eea
\bth
 Let   $A$ be a right $H$-comodule algebra and $M$
     a  right-right (resp.\ left-left) stable \ayd\ module over $H$.
Then  the formulas (\ref{coh1})--(\ref{coh4})
(resp.\  (\ref{hom1})--(\ref{hom4})) define a cyclic module structure on
$\{\hom^H(A^{\ot(n+1)},M)\}_{n\in\N}$ (resp.\ $\{A^{\ot(n+1)}
{\square}_H M\}_{n\in\N}$).
\ethe\noindent
Now one can either define the cyclic cohomology (resp.\ homology) of $A$ with
coefficients in $M$ as the total cohomology (resp.\ homology) of an appropriate
double complex \cite[p.77 and p.72]{l-jl},
or assume that $\Q\inc k$ and proceed as above
Proposition~\ref{pair}.\\

\footnotesize
{\bf Acknowledgements.}
This work was partially supported by
the Marie Curie Fellowship
 HPMF-CT-2000-00523 and KBN grant 2 P03A 013 24 (P.M.H.).
P.M.H.\ is also grateful to the
University of Western Ontario for its hospitality and financial support,
and to the participants of his Hopf-cyclic cohomology workshop for their
interest and encouragement. Special thanks are due to  M.~Furuta and P.~Schauenburg
for their helpful remarks, and to T.~Brzezi\'nski for his help with finding
references.
All  authors are grateful  to R.~Taillefer for the French translation.

 \end{document}